\documentclass[12pt]{amsart}
\usepackage{etex}
\usepackage[english]{babel}
\usepackage{latexsym}
\usepackage{amsbsy}
\usepackage{amssymb}
\usepackage{amsfonts, amsmath}
\usepackage{mathtools}
\usepackage{csquotes}
\usepackage{pst-tree}
\usepackage{pst-math,pst-xkey}
\usepackage{multirow}
\usepackage{blkarray}
\usepackage{tabularx}
\usepackage{arydshln,leftidx,mathtools}
\usepackage{ifthen}
\usepackage{tikz}
\usepackage{tikz-cd}
\usepackage{float}
\usepackage{picins}
\usepackage{amsbsy}
\usetikzlibrary{matrix}
\usetikzlibrary{calc}
\usetikzlibrary{fit}
\usepackage{amsmath,amsfonts,amssymb}
\usepackage{mathrsfs,eurosym}
\usepackage{pst-plot,pst-eucl}
\usepackage{amsfonts}
\usepackage{amsbsy}
\usepackage{alltt}
\usepackage{fontenc}
\usepackage{newlfont}
\usepackage{xlop}
\usepackage{graphicx}
\usepackage{yhmath}
\usepackage{mathdots}

\usepackage{MnSymbol}
\usepackage{amsthm}
\usepackage{url}
\usepackage{enumerate}
\usepackage{leqno}

\newtheorem{thm}{Theorem}[section]
\newtheorem{lem}[thm]{Lemma}

\theoremstyle{definition}

\newtheorem{rmks}[thm]{Remarks}
\newtheorem{ex}[thm]{Example}

\numberwithin{equation}{section}

\usepackage{cite}

\numberwithin{equation}{section}

\DeclarePairedDelimiter\floor{\lfloor}{\rfloor}

\title[Powers of semicirculant and $r$-circulant matrices]{Arbitrary positive powers of semicirculant and $r$-circulant matrices}
\author{\bf M. Mou\c{c}ouf}
\date{}


\begin{document}
\maketitle
\begin{center}
{\footnotesize Department of Mathematics, Faculty of Science, Chouaib Doukkali University, Morocco\\
Email: moucouf@hotmail.com}
\end{center}
\begin{abstract}
We provide a novel recursive method, which does not require any assumption, to compute the entries of the kth power of
a semicirculant matrix. As an application, a method for computing the entries of the kth power of $r$-circulant
matrices is also presented.
\end{abstract}
\section{Introduction}
%
$r$-Circulant matrices comprise an essential topic in numerous areas of mathematics and natural sciences because of their theoretical and applied aspects. Accordingly, they have gained significant importance owing to their frequent use in various applications. They appear in cryptography, number theory, information processing, coding theory, digital image processing, isotropic Markov chain models, spline approximation, among other domains. The positive integer powers for arbitrary circulant and $r$-circulant matrices have been studied~\cite{Feng, Jiang}. This topic was widely studied in literature~\cite{Feng, Jiang, Dav, Guti1, Guti2, Kok, Rim1, Rim2, Rim3, Rim4, Zhao}.
\\An $n\times n$ $r$-circulant matrix $C_{n,r}$ over a unitary commutative ring $R$ has the
following form:
\begin{equation}
C_{n,r}=\left(\begin{array}{cccccc}
c_{0}&c_{1}&c_{2}&\cdots&c_{n-2}&c_{n-1}\\
rc_{n-1}&c_{0}&c_{1}&\cdots&  c_{n-3}   &c_{n-2}\\
rc_{n-2}&rc_{n-1}&c_{0}&\cdots&c_{n-4}&c_{n-3}\\
\cdots&\cdots&\cdots&\cdots&\cdots  &\cdots\\
rc_{1}&rc_{2}&rc_{3}&\cdots&rc_{n-1}&c_{0}
\end{array}
\right).
\end{equation}
The $r$-circulant matrix $C_{n,r}$ is determined by its first-row elements $c_{0},\ldots,c_{n-1}$ and parameter $r\in R$.
Thus, we denote $C_{n,r}=\text{circ}_{n,r}(c_{0},\ldots,c_{n-1})$. When $r=1$, we obtain the
circulant matrix $\text{circ}_{n}(c_{0},\ldots,c_{n-1})$.
\\In~\cite{Feng}, Feng provided a method for computing the kth power of an arbitrary circulant matrix by using the multinomial expansion theorem and considering the fact that a circulant matrix can be expressed as a linear combination of the powers of the basic circulant permutation matrix. Their work was later extended by Jiang~\cite{Jiang} to arbitrary $r$-circulant matrices. The method used by Feng and Jiang is based on the straightforward application of the multinomial expansion theorem and, therefore, requires to solve the difficult problem of determining the set of all $n$-tuples $(k_{0},\ldots,k_{n-1})$ of nonnegative integers that satisfy the following constraints:
\begin{eqnarray}\label{eq;;010}
\begin{cases}
k_{0}+\cdots+k_{n-1}=q\\
k_{1}+2k_{2}+ \cdots + (n-1)k_{n-1}\equiv i\; (\bmod\, n),
\end{cases}
\end{eqnarray}
where $i,q$ and $n$ denote nonnegative integers such that $i\leq n-1$.
\\However, an $n\times n$ matrix $A$ is called semicirculant if it has the following form
(see e.g., Henrici~\cite{Hen} or Davis~\cite{Dav}):
\begin{equation}
A=[a_{0},a_{1},\ldots,a_{n-1}]=
\begin{pmatrix}
a_{0}& a_{1}&. &.&a_{n-1}\\
0&a_{0}&a_{1}&.&a_{n-2}\\
.\;&\ddots &\ddots& \ddots& .\;\;\;\\
.\;& &\ddots& \ddots &a_{1}\\
0&. &. &0&a_{0}
\end{pmatrix}.
\end{equation}
A particular semicirculant matrix is the Jordan block $J=[0,1,0,\ldots,0]$ of size $n$ with the eigenvalue of $0$. A semicirculant matrix $A=[a_{0},a_{1},\ldots,a_{n-1}]$ can be represented as a polynomial in $J$ as follows:
\begin{eqnarray*}
A=a_{0}I_{n}+a_{1}J+\cdots+a_{n-1}J^{n-1},
\end{eqnarray*}
where $I_{n}$ denotes the unit matrix of index n.
\\Let $k$ be any nonnegative integer and $\{a_{n}\}_{n\geq0}$ be a sequence of the elements of a unitary
commutative ring $R$. Accordingly, we have the following~\cite{Wiki}:
\begin{eqnarray*}
[a_{0}, a_{1}, a_{2}, \ldots]^{k}=[a_{0}(k), a_{1}(k), a_{2}(k), \ldots],
\end{eqnarray*}
where
\begin{eqnarray*}
\left\{\begin{array}{lll}
a_{0}(k)&=&a_{0}^{k},\\
a_{m}(k)&=&\displaystyle\frac{1}{ma_{0}}\sum_{i=1}^{m}(ik-m+i)a_{i}a_{m-i}(k),\quad m\geq 1.
\end{array}\right.
\end{eqnarray*}
However, this method is valid only if $m$ and $a_{0}$ are units in ring $R$. Here,
we propose a novel method, which does not require any assumption, to compute the entries of the kth power of a semicirculant matrix.
\\As an application, we present a method for directly calculating the entries of the kth power of an $r$-circulant
matrix $C_{n,r}=\text{circ}_{n,r}(c_{0},\ldots,c_{n-1})$, provided that the entries of the kth power of the associated infinite semicirculant matrix
$[c_{0},\ldots,c_{n-1},0,0,\ldots]$ are known.
\\Precisely, let $r$ be an element of $R$, $n$ be a positive integer, and $T_{n,r}$ be a linear map given as follows:
\begin{eqnarray*}
T_{n,r}(a_{ij})=(r^{\floor*{\frac{|j-i|}{n}}})\odot (a_{ij}),
\end{eqnarray*}
where $\odot$ denotes the Hadamard product, $(a_{ij})$ a finite or an infinite matrix with coefficients
in $R$, and $\floor*{x}$ the greatest integer less than or equal to $x$.
\\Let
 \begin{eqnarray*}
 T_{n,r}([c_{0}, c_{1}, \ldots,c_{n-1},0,0,\ldots]^{k})=[c_{0}(k), c_{1}(k) , c_{2}(k), \ldots].
 \end{eqnarray*}
 Then, one has
\begin{eqnarray*}
C_{n,r}^{k}=\sum_{m\geq 0}\text{circ}_{n,r}(c_{nm}(k),\ldots,c_{n(m+1)-1}(k))
\end{eqnarray*}
In addition to the notation introduced above, we use the notation $\Delta(m,q,p)$ for the solution set of the following system of equations:
\begin{eqnarray*}
\begin{cases}
(k_{1},\ldots,k_{m})\in \mathbb{N}^{m}\\
k_{1}+\cdots+k_{m}=p\\
k_{1}+2k_{2}+ \cdots + mk_{m}=q.
\end{cases}
\end{eqnarray*}
Here, $p,q$ and $m$ denote integers such that
$p\leqslant q\leqslant m$.
If $q=m$, $\Delta(m,q,p)$ will be denoted by $\Delta(m,p)$.
\\We also use notation $[C_{n,r}]$ for the infinite semicirculant matrix $[c_{0}, c_{1}, \ldots,c_{n-1},0,0,\ldots]$
associated with the $r$-circulant matrix $C_{n,r}=\text{circ}_{n,r}(c_{0},\ldots,c_{n-1})$.
\\Throughout this paper, the symbol $R$ will be used to denote an arbitrary commutative ring with identity.
%
\section{Explicit expression of the general terms of a recursive polynomial sequence}
\label{sect:Explicit}
%
Let $R(X_{m})_{m\geq0}$ be an $R$-module
spanned by the family $\{X_{0}, X_{1}, X_{2}, \ldots\}$ of indeterminates over $R$.
\\Let $\nabla : R(X_{m})_{m\geq0}\rightarrow R(X_{m})_{m\geq0}$ be the shift operator defined as
\begin{eqnarray*}
\nabla(X_{m})=X_{m+1}.
\end{eqnarray*}
For a given sequence $\{a_{m}\}_{m\geq0}$ of the elements of a ring $R$, we recursively define the sequence
\begin{eqnarray*}
\{\pmb{a}_{0}(X_{0}), \pmb{a}_{1}(X_{1}), \ldots, \pmb{a}_{m}(X_{1},\ldots,X_{m}), \ldots\}
\end{eqnarray*}
of the elements of $R(X_{m})_{m\geq0}$ as follows:
\begin{equation}\label{eq;5}
\begin{cases}
\pmb{a}_{0}&=X_{0},\\
\pmb{a}_{m+1}&=a_{m+1}a_{0}^{m}\nabla(\pmb{a}_{0})+\cdots+a_{i+1}a_{0}^{i}\nabla(\pmb{a}_{m-i})+\cdots+a_{1}\nabla(\pmb{a}_{m}).
\end{cases}
\end{equation}
It is a routine to verify that the sequence
\begin{eqnarray*}
\{\pmb{a}_{0}(X_{0}), \pmb{a}_{1}(X_{1}), \ldots, \pmb{a}_{m}(X_{1},\ldots,X_{m}), \ldots\}
\end{eqnarray*}
is uniquely determined by the recurrence equation~\eqref{eq;5}.
\\The following result provides an explicit expression for $\pmb{a}_{m}(X_{1},\ldots,X_{m})$.
\begin{thm}\label{Thm;1} The terms of the polynomial sequence defined by the recurrence equation~\eqref{eq;5} assume the following explicit expression:
\begin{equation}\label{eq : 541}
\pmb{a}_{m}=\displaystyle\sum_{p=1}^{m}X_{p}a_{0}^{m-p}\sum_{\substack{\Delta(m,p)}}\binom{p}{k_{1},
\ldots,k_{m}}a_{1}^{k_{1}}\cdots a_{m}^{k_{m}}\quad \text{for}\,\,m\geqslant 1.
\end{equation}
\end{thm}
\proof
For simplicity, let us denote
\begin{eqnarray*}
L(i,p)=\sum_{\substack{\Delta(i,p)}}\binom{p}{k_{1},\ldots,k_{i}}a_{1}^{k_{1}}\cdots a_{i}^{k_{i}}.
\end{eqnarray*}
Substitute
\begin{eqnarray*}
Q_{m}=
\begin{cases}
X_{0} &\quad\text{if}\quad m=0\\
\displaystyle\sum_{p=1}^{m}X_{p}a_{0}^{m-p}\sum_{\substack{\Delta(m,p)}}\binom{p}{k_{1},\ldots,k_{m}}a_{1}^{k_{1}}
\cdots a_{m}^{k_{m}} &\quad\text{otherwise}.
\end{cases}
\end{eqnarray*}
We must prove that sequences $\{\pmb{a}_{m}\}_{m}$ and $\{Q_{m}\}_{m}$ are identical. Because sequence $\{\pmb{a}_{m}\}_{m}$ is uniquely determined by~\eqref{eq;5} and $\{Q_{m}\}_{m}$ begins with $X_{0}$, it suffices to show that
\begin{eqnarray*}
Q_{m+1}=a_{m+1}a_{0}^{m}X_{1}+\displaystyle\sum_{p=1}^{m}a_{0}^{m-p}a_{m-p+1}\nabla(Q_{p}).
\end{eqnarray*}
Accordingly, we argue as follows. For any positive integer $p\leq m$, we have the following:
\begin{equation}\label{eq.01}
Q_{p}=\displaystyle\sum_{i=1}^{p}X_{i}a_{0}^{p-i}L(p,i),
\end{equation}
and
\begin{eqnarray*}
a_{0}^{m-p}a_{m-p+1}\nabla(Q_{p})=\displaystyle\sum_{i=1}^{p}X_{i+1}a_{0}^{m-i}L(p,i)a_{m-p+1}.
\end{eqnarray*}
Thus, for all nonnegative integers $m$, we have
\begin{eqnarray*}
a_{m+1}a_{0}^{m}X_{1}+\displaystyle\sum_{p=1}^{m}a_{0}^{m-p}a_{m-p+1}\nabla(Q_{p})=
a_{m+1}a_{0}^{m}X_{1}+\displaystyle\sum_{p=1}^{m}\sum_{i=1}^{p}X_{i+1}a_{0}^{m-i}L(p,i)a_{m-p+1}.
\end{eqnarray*}
However, it is clear that
\begin{eqnarray*}
\begin{array}{l}
\displaystyle\sum_{p=1}^{m}\sum_{i=1}^{p}X_{i+1}a_{0}^{m-i}L(p,i)a_{m-p+1}=\\
\hspace*{2cm}\begin{array}{ccccccc}
X_{2}a_{0}^{m-1}L(1,1)a_{m}&+&&&&\\
X_{2}a_{0}^{m-1}L(2,1)a_{m-1}&+&X_{3}a_{0}^{m-2}L(2,2)a_{m-1}&+& &&\\
\vdots                       & & \vdots                      & & &&\\
X_{2}a_{0}^{m-1}L(m,1)a_{1}&+&X_{3}a_{0}^{m-2}L(m,2)a_{1}&+&\cdots&+&X_{m+1}L(m,m)a_{1}.
\end{array}
\end{array}
\end{eqnarray*}
By rearranging the terms in the right-hand side of this equation, one obtains
\begin{equation}\label{eq.02}
\displaystyle\sum_{p=1}^{m}\sum_{i=1}^{p}X_{i+1}a_{0}^{m-i}L(p,i)a_{m-p+1}=
\displaystyle\sum_{p=1}^{m}X_{p+1}a_{0}^{m-p}\sum_{i=p}^{m}L(i,p)a_{m-i+1}.
\end{equation}
Because $L(m+1,1)=a_{m+1}$, we have from~\eqref{eq.01} the following:
\begin{eqnarray}\label{eq.03}
Q_{m+1}&=&a_{m+1}a_{0}^{m}X_{1}+\displaystyle\sum_{i=2}^{m+1}X_{i}a_{0}^{m+1-i}L(m+1,i)\nonumber\\
&=&a_{m+1}a_{0}^{m}X_{1}+\displaystyle\sum_{p=1}^{m}X_{p+1}a_{0}^{m-p}L(m+1,p+1).
\end{eqnarray}
Comparing \eqref{eq.02} with \eqref{eq.03}, it is evident that we only need to prove the following:
\begin{eqnarray*}
L(m+1,p+1)=\sum_{i=p}^{m}L(i,p)a_{m-i+1}.
\end{eqnarray*}
Thus, the proof will be complete with the lemma below.
\endproof
\begin{lem} Let $\{a_{n}\}_{n\geq1}$ be a sequence of the elements of $R$, and
$\{L(i,p)\}_{1\leqslant p\leqslant i}$ be a double sequence defined as
\begin{eqnarray*}
  L(i,p)=\sum_{\substack{\Delta(i,p)}}\binom{p}{k_{1},\ldots,k_{i}}a_{1}^{k_{1}}\cdots a_{i}^{k_{i}}.
\end{eqnarray*}
Accordingly, for all $m\geqslant p\geqslant 1$, one has
\begin{eqnarray*}
L(m+1,p+1)=\sum_{i=1}^{m-p+1}L(m-i+1,p)a_{i}.
\end{eqnarray*}
\end{lem}
\proof Consider the following system:
\begin{eqnarray*}
\begin{cases}
(k_{1},\ldots,k_{m+1})\in \mathbb{N}^{m+1}\\
k_{1}+\cdots+k_{m+1}=p \\
k_{1}+2k_{2}+\cdots+(m+1)k_{m+1}=q .
\end{cases}
\end{eqnarray*}
 By subtracting the two equations of the above system, we obtain
 \begin{eqnarray*}
 k_{2}+\cdots+(q-p+1)k_{q-p+1}+\cdots+mk_{m+1}=q-p.
 \end{eqnarray*}
 It follows that $k_{j}=0$ for all $j\geqslant q-p+1$, and hence $k_{j}=0$ for all $j\geqslant q+1$. Thus, the following
equivalence is valid:
\begin{eqnarray}\label{eq;;1}
(k_{1},\ldots,k_{m+1})\in \Delta(m+1,q,p)\Longleftrightarrow
\begin{cases}
(k_{1},\ldots,k_{q})\in \Delta(q,p)\\
 k_{q+1}=\cdots=k_{m+1}=0.
 \end{cases}
\end{eqnarray}
Now consider polynomial $Q(X)=a_{1}X+\cdots+a_{n}X^{n}$.
\\Using the multinomial theorem we obtain
\begin{eqnarray*}
Q(X)^{p}&=&(a_{1}X+\cdots+a_{n}X^{n})^{p}\\
&=&\displaystyle\sum_{k_{1}+\cdots+k_{n}=p}
\binom{p}{k_{1},\ldots,k_{n}}a_{1}^{k_{1}}\cdots a_{n}^{k_{n}}X^{k_{1}+2k_{2}+\cdots+nk_{n}}\\
&=&\sum_{i=p}^{np} \sum_{\Delta(n,i,p)} \binom{p}{k_{1},\ldots,k_{n}}a_{1}^{k_{1}}\cdots a_{n}^{k_{n}}X^{i},
\end{eqnarray*}
and considering~\eqref{eq;;1}, we obtain
\begin{eqnarray*}
Q(X)^{p}&=&\sum_{i=p}^{np} \sum_{\Delta(i,p)} \binom{p}{k_{1},\ldots,k_{i}}a_{1}^{k_{1}}\cdots a_{i}^{k_{i}}X^{i}\\
&=&\sum_{i=p}^{np} L(i,p)X^{i}.
\end{eqnarray*}
We can then write $Q(X)^{p+1}$ as
\begin{eqnarray*}
Q(X)^{p+1}=\sum_{i=p+1}^{n(p+1)} L(i,p+1)X^{i}.
\end{eqnarray*}
Now, if we rewrite $Q(X)^{p+1}$ as
\begin{eqnarray*}
Q(X)^{p+1}=Q(X)^{p}(a_{1}X+\cdots+a_{n}X^{n}),
\end{eqnarray*}
 the term of degree $m+1$ of $Q(X)^{p+1}$
is the same as that of the polynomial \begin{eqnarray*}
(\sum_{i=1}^{m-i+1}L(m-i+1,p)X^{m-i+1})(a_{1}X+\cdots+a_{n}X^{n}),
\end{eqnarray*}
that is
\begin{eqnarray*}
L(m+1,p+1)X^{m+1}&=&\sum_{i=1}^{m-p+1}L(m-i+1,p)X^{m-i+1}a_{i}X^{i}\\
&=&\sum_{i=1}^{m-p+1}L(m-i+1,p)a_{i}X^{m+1}.
\end{eqnarray*}
Therefore,
\begin{eqnarray*}
L(m+1,p+1)=\sum_{i=1}^{m-p+1}L(m-i+1,p)a_{i},
\end{eqnarray*}
and the lemma follows.
\endproof
%
\section{Simple recursive formula for computing the kth power of semicirculant matrices}
\label{sect:A recursive formula}
%
From Formula~\ref{eq : 541}, we have
\begin{eqnarray*}
\pmb{a}_{m}(X_{1},\ldots,X_{m})=\sum_{p=1}^{m}a_{0}^{m-p}L(m,p)X_{p},
\end{eqnarray*}
where
\begin{eqnarray*}
L(m,p)=\sum_{\substack{\Delta(m,p)}}\binom{p}{k_{1},\ldots,k_{m}}a_{1}^{k_{1}}\cdots a_{m}^{k_{m}}.
\end{eqnarray*}
Next, we denote by $\{\pmb{a}_{m}^{(k)}\}_{m,k\geq0}$ the double sequence over $R$ defined as follows:
\begin{equation}\label{eq;1}
\pmb{a}_{m}^{(k)}=\begin{cases}
a_{0}^{k}\binom{k}{0} \quad\text{if}\quad m=0\\
\sum_{p=1}^{m}L(m,p)a_{0}^{k-p}\binom{k}{p} \quad\text{otherwise}.
\end{cases}
\end{equation}
Clearly, if $a_{0}$ is a unit in $R$, then
\begin{eqnarray*}
\pmb{a}_{m}^{(k)}=a_{0}^{k-m}\pmb{a}_{m}(\binom{k}{1},\ldots,\binom{k}{m}),
\end{eqnarray*}
for all nonnegative integers $k$.
\\Hereinafter, we adopt the convention that for any element $a\in R$
and any nonnegative integers $k\leq p$,
\begin{equation}\label{eq,,1}
a^{k-p}\binom{k}{p}=\delta_{k,p},
\end{equation}
where
\begin{eqnarray*}
\delta_{k,p}=\begin{cases}
1 \quad\text{if}\quad k=p\\0 \quad\text{otherwise}
\end{cases}
\end{eqnarray*}
\\Using the above-mentioned notation, we state the following result.
\begin{thm} Let $\{a_{n}\}_{n\geq0}$ be a sequence of the elements of a ring $R$.
For all nonnegative integers $k$, the entries of the kth power of the infinite semicirculant matrix
$A=[a_{0}, a_{1} ,a_{2}, \ldots]$ are given as follows:
\begin{eqnarray*}
A^{k}=[\pmb{a}_{0}^{(k)}, \pmb{a}_{1}^{(k)}, \pmb{a}_{2}^{(k)}, \ldots].
\end{eqnarray*}
\end{thm}
\proof From formula~(1.3-3) of~\cite{Hen}, we may assume, without loss of generality, that
$A=[a_{0},a_{1},\ldots,a_{n}]$ is a finite semicirculant matrix. Let $A^{k}=[a_{0}(k),a_{1}(k),\ldots,a_{n}(k)]$.
We wish to prove that for any nonnegative integer $m\leq n$, $a_{m}(k)=\pmb{a}_{m}^{(k)}$.
\\Because the claim is true for $m=0$, assume that $m\geqslant 1$.
Because
\begin{eqnarray*}A^{k}&=&(a_{0}I_{n+1}+a_{1}J+\cdots+a_{n}J^{n})^{k}\\
&=&\displaystyle\sum_{\substack{k_{0},\ldots,k_{n}\geqslant 0\\k_{0}+\cdots+k_{n}=k}}
\binom{k}{k_{0},\ldots,k_{n}}a_{0}^{k_{0}}a_{1}^{k_{1}}\cdots a_{n}^{k_{n}}J^{k_{1}+2k_{2}+\cdots+nk_{n}},
\end{eqnarray*}
we have
\begin{eqnarray*}
a_{m}(k)&=&\displaystyle\sum_{(k_{1},\ldots,k_{n})\in\substack{\Delta(n,m,k-k_{0})}}
\binom{k}{k_{0},\ldots,k_{n}}a_{0}^{k_{0}}a_{1}^{k_{1}}\cdots a_{n}^{k_{n}}\\
 &=&\displaystyle\sum_{(k_{1},\ldots,k_{m})\in\substack{\Delta(m,k-k_{0})}}
\binom{k}{k_{0},\ldots,k_{m}}a_{0}^{k_{0}}a_{1}^{k_{1}}\cdots a_{m}^{k_{m}}.\\
\end{eqnarray*}
Therefore,
\begin{eqnarray*}
a_{m}(k)&=&\sum_{p=1}^{m}a_{0}^{k-p}\binom{k}{p}\sum_{\substack{\Delta(m,p)}}
\binom{p}{k_{1},\ldots,k_{m}}a_{1}^{k_{1}}\cdots a_{m}^{k_{m}}\\
&=&\pmb{a}_{m}^{(k)}.
\end{eqnarray*}
Thus, the proof is completed.
\endproof
Next, we provide a recursive method to directly compute the coefficients $a_{i}(k)$ of $[a_{0}, a_{1}, a_{2}, \ldots]^{k}$
without using polynomials $\pmb{a}_{i}$; however, only the terms $\binom{k}{p}$ are considered indeterminates.
\\To do this, let us transform the double sequence $\{\pmb{a}_{m}^{(k)}\}_{m,k}$, which is defined in \eqref{eq;1}, into another
double sequence $\{\pmb{\widehat{a}}_{m}^{(k)}\}_{m,k}$, which we define as follows:
\begin{equation}
\pmb{\widehat{a}}_{m}^{(k)}=
\begin{cases}
a_{0}^{k-1}\binom{k}{1} \quad\text{if}\quad m=0\\
\sum_{p=1}^{m}L(m,p)a_{0}^{k-p-1}\binom{k}{p+1} \quad\text{otherwise}.
\end{cases}
\end{equation}
The term $\pmb{\widehat{a}}_{m}^{(k)}$ can be obtained from $\pmb{a}_{m}^{(k)}$ by replacing the integer $p$ in each sequence $a_{0}^{k-p}\binom{k}{p}$ by $p+1$.
\\Next, we introduce the following functions of variable $x$:
\begin{eqnarray*}
g_{m,k}(x)=
\begin{cases}
x^{k-1}\binom{k}{1} \quad\text{if}\quad m=0\\
\sum_{p=1}^{m}L(m,p)x^{k-p-1}\binom{k}{p+1} \quad\text{otherwise},\\
x\in R
\end{cases}
\end{eqnarray*}
and
\begin{eqnarray*}
f_{m,k}(x)=
\begin{cases}
x^{k}\binom{k}{0} \quad\text{if}\quad m=0\\
x^{k-m}\sum_{p=1}^{m}L(m,p)x^{m-p}\binom{k}{p} \quad\text{otherwise}.\\
x\in R
\end{cases}
\end{eqnarray*}
Because $f_{m,k}(x)=\sum_{p=1}^{m}L(m,p)x^{k-p}\binom{k}{p}$, according to the convention adopted above~\eqref{eq,,1}, functions $f_{m,k}$ are defined for all $x\in R$. Therefore, we write $\pmb{a}_{m}^{(k)}$ in terms
of $\pmb{a}_{m}(\binom{k}{1},\ldots,\binom{k}{m})=\sum_{p=1}^{m}L(m,p)a_{0}^{m-p}\binom{k}{p}$ even when $a_{0}$ is not a unit in ring $R$, as follows: $\pmb{a}_{m}^{(k)}=f_{m,k}(a_{0})$.
\\Next, we consider the sequence $\{x_{m}\}_{m\geq0}$ defined by $x_{0}=x$ and $x_{m}=a_{m}$ for all $m\geq 1$.
From Theorem~\ref{Thm;1}, one has
\begin{eqnarray*}
f_{m,k}(x)=
\begin{cases}
x^{k}\pmb{x}_{0}(\binom{k}{0}) \quad\text{if}\quad m=0\\
x^{k-m}\pmb{x}_{m}(\binom{k}{1},\ldots,\binom{k}{m}) \quad\text{otherwise}
\end{cases}
\end{eqnarray*}
and
\begin{eqnarray*}
g_{m,k}(x)=
\begin{cases}
x^{k-1}\nabla(\pmb{x}_{0})(\binom{k}{1}) \quad\text{if}\quad m=0\\
x^{k-m-1}\nabla(\pmb{x}_{m})(\binom{k}{2},\ldots,\binom{k}{m+1}) \quad\text{otherwise}.
\end{cases}
\end{eqnarray*}
Thus,
\begin{eqnarray*}
f_{m+1,k}(x)&=&x^{k-m-1}\pmb{x}_{m+1}(\binom{k}{1},\ldots,\binom{k}{m+1})\\
&=&x^{k-m-1}[a_{m+1}x_{0}^{m}\nabla(\pmb{x}_{0})(\binom{k}{1})+\cdots+
a_{i+1}x_{0}^{i}\nabla(\pmb{x}_{m-i})(\binom{k}{1},\ldots,\binom{k}{m-i+1})+
\cdots\\
&&+a_{1}\nabla(\pmb{x}_{m})(\binom{k}{1},\ldots,\binom{k}{m+1})]\\
&=&a_{m+1}x_{0}^{k-1}\nabla(\pmb{x}_{0})(\binom{k}{1})+\cdots+
a_{i+1}x_{0}^{k-m+i-1}\nabla(\pmb{x}_{m-i})(\binom{k}{1},\ldots,\binom{k}{m-i+1})+\cdots\\
&&+a_{1}x^{k-m-1}\nabla(\pmb{x}_{m})(\binom{k}{1},\ldots,\binom{k}{m+1}).
\end{eqnarray*}
Therefore,
\begin{equation}\label{eq;25}
f_{m+1,k}(x)=a_{m+1}g_{m,k}(x)+\cdots+a_{1}g_{0,k}(x).
\end{equation}
Considering $x=a_{0}$ in~\eqref{eq;25}, we obtain the following recursive formula:
\begin{equation}
\pmb{a}_{m+1}^{(k)}=a_{m+1}\pmb{\widehat{a}}_{0}^{(k)}\cdots+a_{i+1}\pmb{\widehat{a}}_{m-i}^{(k)}+
\cdots+a_{1}\pmb{\widehat{a}}_{m}^{(k)}.
\end{equation}
Conclusively, we have proved the following theorem.
\begin{thm}\label{th3} The double sequence $\{\pmb{a}_{n}^{(k)}\}_{n,k\geq0}$ defined in~\eqref{eq;1} can be recursively determined
as follows:
\begin{equation}\label{equ 2}
\pmb{a}_{m+1}^{(k)}=a_{m+1}\pmb{\widehat{a}}_{0}^{(k)}\cdots+a_{i+1}\pmb{\widehat{a}}_{m-i}^{(k)}+
\cdots+a_{1}\pmb{\widehat{a}}_{m}^{(k)}
\end{equation}
where $\pmb{\widehat{a}}_{n}^{(k)}$ denotes the transform of the term $\pmb{a}_{n}^{(k)}$ obtained by replacing the integer $p$ in each sequence $a_{0}^{k-p}\binom{k}{p}$ by $p+1$.
\end{thm}
\begin{rmks}~
\begin{enumerate}[1.]
\item Notably, $[\pmb{a}_{1}^{(k)}, \pmb{a}_{2}^{(k)}, \ldots]=
[a_{1}, a_{2}, \ldots][\pmb{\widehat{a}}_{0}^{(k)}, \pmb{\widehat{a}}_{1}^{(k)}, \ldots]$. This formula can
help us remember the manner to compute the double sequence $\pmb{a}_{m}^{(k)}$.
\item Notably, in formula~$\eqref{equ 2}$, the elements
$\binom{k}{m}$, $m\geq 0$, must be considered independent indeterminates over ring $R$ during all the formal operations,
and only restored as binomial coefficients when operations are completed.
\item Notably, for some nonnegative integer $m$, $(\pmb{a}_{m}^{(k)})_{k\geq0}$ is an identically zero sequence, then
polynomial $\pmb{a}_{m}$ is identically equal to $0$ and, therefore, so is its transform $\nabla(\pmb{a}_{m})$.
Particularly,
$(\pmb{\widehat{a}}_{m}^{(k)})_{k\geq0}$ is the identically zero sequence.
\item If matrix $C$ has the form $[0,\ldots,0,a_{p},a_{p+1},\ldots]$, it becomes significantly easy to begin by computing the powers of the left-shift matrix of $C$ by $p$ positions and then deduce those of matrix $C$.
\item Clearly, if $a_{n}=0$ for all $n\geq p+1$, then
$\pmb{a}_{m}^{(k)}=\displaystyle\sum_{i=\lceil\frac{m}{p}\rceil}^{m} L(m,i)a_{0}^{k-i}\binom{k}{i}$, where $\lceil x\rceil$
denotes the smallest integer greater than or equal to $x$. Particularly,
$\pmb{a}_{m}^{(k)}=0$ for all $m\geq kp+1$.
\end{enumerate}
\end{rmks}
The following two simple examples illustrate the application of Theorem~\ref{th3}.
\begin{ex}~
\\Consider a semicirculant matrix $A=[2,4,2,3]$ over the ring $\mathbb{Z}/8\mathbb{Z}$ of integers modulo $8$. Let $k$ be any nonnegative integer. Accordingly, $A^{k}=[\pmb{a}_{0}^{(k)}, \pmb{a}_{1}^{(k)}, \pmb{a}_{2}^{(k)}, \pmb{a}_{3}^{(k)}]$, where
\begin{eqnarray*}
\pmb{a}_{0}^{(k)}&=&2^{k}\binom{k}{0}\\
\pmb{\widehat{a}}_{0}^{(k)}&=&2^{k-1}\binom{k}{1};\\
\pmb{a}_{1}^{(k)}&=&4\times2^{k-1}\binom{k}{1}=2^{k+1}\binom{k}{1}\\
\pmb{\widehat{a}}_{1}^{(k)}&=&2^{k}\binom{k}{2};\\
\pmb{a}_{2}^{(k)}&=&4\times2^{k}\binom{k}{2}+2\times2^{k-1}\binom{k}{1}=2^{k+2}\binom{k}{2}+2^{k}\binom{k}{1}\\
\pmb{\widehat{a}}_{2}^{(k)}&=&2^{k+1}\binom{k}{3}+2^{k-1}\binom{k}{2};\\
\pmb{a}_{3}^{(k)}&=&4\times(2^{k+1}\binom{k}{3}+2^{k-1}\binom{k}{2})+2\times2^{k}\binom{k}{2}+3\times2^{k-1}\binom{k}{1}\\
&=&2^{k+3}\binom{k}{3}+2^{k+2}\binom{k}{2}+3\times2^{k-1}\binom{k}{1}
\end{eqnarray*}
\end{ex}
\begin{ex}~
\\Consider another semicirculant matrix $B=[0,2,1,1,0]$ over the ring $\mathbb{Z}/8\mathbb{Z}$ of integers modulo $8$. Let $k$ be any nonnegative integer. Accordingly, $B^{k}=[\pmb{b}_{0}^{(k)}, \pmb{b}_{1}^{(k)}, \pmb{b}_{2}^{(k)}, \pmb{b}_{3}^{(k)}]$, where
\begin{eqnarray*}
\pmb{b}_{0}^{(k)}&=&0^{k}\binom{k}{0}\\
\pmb{\widehat{b}}_{0}^{(k)}&=&0^{k-1}\binom{k}{1};\\
\pmb{b}_{1}^{(k)}&=&2\times0^{k-1}\binom{k}{1}\\
\pmb{\widehat{b}}_{1}^{(k)}&=&2\times0^{k-2}\binom{k}{2};\\
\pmb{b}_{2}^{(k)}&=&4\times0^{k-2}\binom{k}{2}+0^{k-1}\binom{k}{1}\\
\pmb{\widehat{b}}_{2}^{(k)}&=&4\times0^{k-3}\binom{k}{3}+0^{k-2}\binom{k}{2};\\
\pmb{b}_{3}^{(k)}&=&2\times(4\times0^{k-3}\binom{k}{3}+0^{k-2}\binom{k}{2})+2\times0^{k-2}\binom{k}{2}+0^{k-1}\binom{k}{1}
\end{eqnarray*}
\begin{eqnarray*}
&=&8\times0^{k-3}\binom{k}{3}+4\times0^{k-2}\binom{k}{2}+0^{k-1}\binom{k}{1}\\
\pmb{\widehat{b}}_{3}^{(k)}&=&8\times0^{k-4}\binom{k}{4}+4\times0^{k-3}\binom{k}{3}+0^{k-2}\binom{k}{2};\\
\pmb{b}_{4}^{(k)}&=&2\times(8\times0^{k-4}\binom{k}{4}+4\times0^{k-3}\binom{k}{3}+0^{k-2}\binom{k}{2})+
4\times0^{k-3}\binom{k}{3}+\\
&&0^{k-2}\binom{k}{2}+2\times0^{k-2}\binom{k}{2}+0\times0^{k-1}\binom{k}{1}\\
&=&16\times0^{k-4}\binom{k}{4}+12\times0^{k-3}\binom{k}{3}+5\times0^{k-2}\binom{k}{2}.
\end{eqnarray*}
\end{ex}
%
\section{ kth power of matrix $C_{n,r}$}
\label{sect:The calculation}
%
Let $E_{n,r}$ be the basic $r$-circulant permutation matrix of order $n$ over ring $R$ defined as
$E_{n,r}=\text{circ}_{n,r}(0,1,0,\ldots,0)$. Because $E_{n,r}$ is the companion matrix of polynomial
$X^{n}-r$, it follows that $X^{n}-r$ is the characteristic polynomial of matrix  $E_{n,r}$. Accordingly, we have
$E_{n,r}^{n}=r$ (see, e.g., Brown~\cite{Brow}). Using this result, we may deduce that $E_{n,r}^{k+n}=rE_{n,r}^{k}$
for all $k\geqslant 0$. Hence,
\begin{equation}\label{eq 1111}
E_{n,r}^{k}=\left\{\begin{array}{ccc}
r^{\floor*{\frac{k}{n}}}E_{n,r}^{s},&\text{if}& k\equiv s \; (\bmod\, n)\\
r^{\floor*{\frac{k}{n}}}I_{n},&\text{if}& k\equiv 0 \; (\bmod\, n)
\end{array}\right.
\end{equation}
where $1\leq r\leq n-1$.
\\The next theorem shows that the entries of the kth power of any $r$-circulant matrix $C_{n,r}$ can be deduced from those
of the kth power of the infinite semicirculant matrix $[C_{n,r}]$.
\begin{thm}\label{Thm C} Let $C_{n,r}=\text{circ}_{n,r}(c_{0},\ldots,c_{n-1})$ be an $r$-circulant matrix and let
\begin{eqnarray*}
[C_{n,r}]=[c_{0},\ldots,c_{n-1},0,0,\ldots]
\end{eqnarray*}
be the associated infinite semicirculant matrix. Put
\begin{eqnarray*}
 T_{n,r}([C_{n,r}]^{k})=[b_{0}(k), b_{1}(k), b_{2}(k), \ldots],
\end{eqnarray*}
 where $k$ denotes any nonnegative integer and $T_{n,r}$ a linear map defined as
\begin{eqnarray*}
T_{n,r}(a_{ij})=(r^{\floor*{\frac{|j-i|}{n}}})\odot (a_{ij}).
\end{eqnarray*}
  Accordingly, we have
\begin{eqnarray*}
C_{n,r}^{k}=\sum_{m\geq 0}\text{circ}_{n,r}(c_{nm}(k),\ldots,c_{n(m+1)-1}(k)).
\end{eqnarray*}
\end{thm}
\proof  It is well-known that
\begin{eqnarray*}
C_{n,r}=c_{0}I_{n}+c_{1}E_{n,r}+ \cdots +c_{n-1}E_{n,r}^{n-1}.
\end{eqnarray*}
Therefore, for any nonnegative integer $k$,
\begin{eqnarray*}
C_{n,r}^{k}=(c_{0}I_{n}+c_{1}E_{n,r}+ \cdots +c_{n-1}E_{n,r}^{n-1})^{k}.
\end{eqnarray*}
Using the multinomial theorem, we have
\begin{eqnarray*}
C_{n,r}^{k}=\sum_{\substack{k_{0},\ldots,k_{n-1}\in \mathbb{N}\\k_{0}+\cdots+k_{n-1}=k}}
\binom{k}{k_{0},\ldots,k_{n-1}}c_{0}^{k_{0}}c_{1}^{k_{1}}\cdots c_{n-1}^{k_{n-1}}E_{n,r}^{k_{1}+2k_{2}+\cdots+(n-1)k_{n-1}}.
\end{eqnarray*}
Because the family
\begin{eqnarray*}
\{(k_{0},\ldots,k_{n-1})\in \mathbb{N}^{n}/ \,k_{0}+\cdots+k_{n-1}=k \,\,\text{and}\,\, k_{1}+ \cdots +(n-1)k_{n-1}=i+mn\},
m=0,1,\ldots
\end{eqnarray*}
forms a partition of the set
\begin{eqnarray*}
\{(k_{0},\ldots,k_{n-1})\in \mathbb{N}^{n}/ \,k_{0}+\cdots+k_{n-1}=k \,\,\text{and}\,\,
k_{1}+\cdots+(n-1)k_{n-1}\,\,\equiv\,\, i \; (\bmod\, n)\},
\end{eqnarray*}
it follows from~\eqref{eq 1111} that
\begin{eqnarray*}
C_{n,r}^{k}=\sum_{p=0}^{n-1}\sum_{\substack{k_{0}+\cdots+k_{n-1}=k\\
m=k_{1}+\cdots+(n-1)k_{n-1} \,\,\equiv\,\, p \; (\bmod\, n)}}\binom{k}{k_{0},\ldots,k_{n-1}}c_{0}^{k_{0}}c_{1}^{k_{1}}
\cdots c_{n-1}^{k_{n-1}}r^{\floor*{\frac{m}{n}}}E_{n,r}^{p}.
\end{eqnarray*}
Consequently, for all $p=0,\ldots,n-1$, we have
\begin{eqnarray*}
(C_{n,r}^{k})_{p}=\sum_{\substack{k_{0}+\cdots+k_{n-1}=k\\m=k_{1}+\cdots+(n-1)k_{n-1} \,\,\equiv\,\, p \; (\bmod\, n)}}
\binom{k}{k_{0},\ldots,k_{n-1}}c_{0}^{k_{0}}c_{1}^{k_{1}}\cdots c_{n-1}^{k_{n-1}}r^{\floor*{\frac{m}{n}}},
\end{eqnarray*}
where $(C_{n,r}^{k})_{p}$ denotes the pth strip of matrix $C_{n,r}^{k}$.
The conclusion immediately follows from the fact that \begin{eqnarray*}
\sum_{\substack{k_{0}+\cdots+k_{n-1}=k\\k_{1}+\cdots+(n-1)k_{n-1}=m}}
\binom{k}{k_{0},\ldots,k_{n-1}}c_{0}^{k_{0}}c_{1}^{k_{1}}\cdots c_{n-1}^{k_{n-1}}
\end{eqnarray*}
is the $(1,m+1)$ entry of matrix $[C_{n,r}]^{k}$.
\endproof
\begin{rmks}~
\begin{enumerate}[1.]
\item The canonical mapping
\begin{eqnarray*}
c_{0}+c_{1}x+c_{2}x^{2}+\cdots\longrightarrow [c_{0}, c_{1}, c_{2}, \ldots]
\end{eqnarray*}
is an isomorphism from the ring of formal power series onto the ring of infinite semicirculant matrices. Thus, another way of establishing Theorem~\ref{Thm C} is to observe that $C_{n,r}^{k}=P^{k}(E_{n,r})$, where $P(x)$ denotes the representer polynomial of the $r$-circulant matrix $C_{n,r}$.
\item The sequence $r^{\floor*{\frac{m}{n}}}$, which appears in the pth strip
\begin{eqnarray*}
(C_{n,r}^{k})_{p}=\sum_{m\equiv p}c_{m}(k)r^{\floor*{\frac{m}{n}}}
\end{eqnarray*}
of the $r$-circulant matrix $C_{n,r}^{k}$, is a geometric sequence with
the common ratio of $r$.
\end{enumerate}
\end{rmks}
Let us illustrate our method using the following examples. The first example is similar to the one used in~\cite{Jiang}.
\begin{ex} Let $C=circ_{5,-1}(5, 4, 3, 2, 1)$ and let
$[C]=[5, 4, 3, 2, 1,0,0,\ldots]$ be the associated infinite semicirculant matrix. Put
\begin{eqnarray*}
[C]^{k}=[\pmb{c}_{0}^{(k)}, \pmb{c}_{1}^{(k)}, \pmb{c}_{2}^{(k)}, \ldots].
\end{eqnarray*}
Using the method provided in Theorem~\ref{th3}, we obtain
\begin{eqnarray*}
\pmb{c}_{0}^{(k)}&=&5^{k}\binom{k}{0}; \quad\pmb{c}_{0}^{(3)}=125\\
\pmb{c}_{1}^{(k)}&=&4\times5^{k-1}\binom{k}{1}; \quad\pmb{c}_{1}^{(3)}=300\\
\pmb{c}_{2}^{(k)}&=&4^{2}\times5^{k-2}\binom{k}{2}+3\times5^{k-1}\binom{k}{1}; \quad\pmb{c}_{2}^{(3)}=465\\
\pmb{c}_{3}^{(k)}&=&4(4^{2}\times5^{k-3}\binom{k}{3}+3\times5^{k-2}\binom{k}{2})+3\times4\times5^{k-2}\binom{k}{2}+\\
&&2\times5^{k-1}\binom{k}{1}; \quad\pmb{c}_{3}^{(3)}=574\\
\pmb{c}_{4}^{(k)}&=&4\widehat{c}_{3}(k)+3\widehat{c}_{2}(k)+2\widehat{c}_{1}(k)+\widehat{c}_{0}(k); \quad\pmb{c}_{4}^{(3)}=594\\
\pmb{c}_{5}^{(k)}&=&4\widehat{c}_{4}(k)+3\widehat{c}_{3}(k)+2\widehat{c}_{2}(k)+\widehat{c}_{1}(k); \pmb{c}_{5}^{(3)}=504\\
\pmb{c}_{6}^{(k)}&=&4\widehat{c}_{5}(k)+3\widehat{c}_{4}(k)+2\widehat{c}_{3}(k)+\widehat{c}_{2}(k); \quad\pmb{c}_{6}^{(3)}=369\\
\pmb{c}_{7}^{(k)}&=&4\widehat{c}_{6}(k)+3\widehat{c}_{5}(k)+2\widehat{c}_{4}(k)+\widehat{c}_{3}(k); \quad\pmb{c}_{7}^{(3)}=234\\
\pmb{c}_{8}^{(k)}&=&4\widehat{c}_{7}(k)+3\widehat{c}_{6}(k)+2\widehat{c}_{5}(k)+\widehat{c}_{4}(k); \quad\pmb{c}_{8}^{(3)}=126\\
\pmb{c}_{9}^{(k)}&=&4\widehat{c}_{8}(k)+3\widehat{c}_{7}(k)+2\widehat{c}_{6}(k)+\widehat{c}_{5}(k); \quad\pmb{c}_{9}^{(3)}=56\\
\pmb{c}_{10}^{(k)}&=&4\widehat{c}_{9}(k)+3\widehat{c}_{8}(k)+2\widehat{c}_{7}(k)+\widehat{c}_{6}(k);
\quad\pmb{c}_{10}^{(3)}=21\\
\pmb{c}_{11}^{(k)}&=&4\widehat{c}_{10}(k)+3\widehat{c}_{9}(k)+2\widehat{c}_{8}(k)+\widehat{c}_{7}(k);
\quad\pmb{c}_{11}^{(3)}=6\\
\pmb{c}_{12}^{(k)}&=&4\widehat{c}_{11}(k)+3\widehat{c}_{10}(k)+2\widehat{c}_{9}(k)+\widehat{c}_{8}(k); \quad\pmb{c}_{12}^{(3)}=1\\
\pmb{c}_{s}^{(3)}&=&0 \quad\text{for all}\quad s\geq (3\times4)+1=13.
\end{eqnarray*}
Hence, \\$[C]^{3}=[125, 300, 465, 574, 594, 504, 369, 234, 126, 56, 21, 6, 1,0,0,\ldots]$.
\\Therefore,
\begin{eqnarray*}
C_{0}^{3}&=&125(-1)^{0}+504(-1)^{1}+21(-1)^{2}=-358\\
C_{1}^{3}&=&300(-1)^{0}+369(-1)^{1}+6(-1)^{2}=-63\\
C_{2}^{3}&=&465(-1)^{0}+234(-1)^{1}+(-1)^{2}=232\\
C_{3}^{3}&=&574(-1)^{0}+126(-1)^{1}=448\\
C_{4}^{3}&=&594(-1)^{0}+56(-1)^{1}=538.
\end{eqnarray*}
Thus,
\begin{eqnarray*}
C^{3}=circ_{5,-1}(-358, -63, 232, 448, 538).
\end{eqnarray*}
\end{ex}
\begin{ex} Let $p,q$ and $n$ denote integers such that $0\leq p<q<n$ and let $a,b\in R$. Consider an $r$-circulant matrix
$C_{n,r}=\text{circ}_{n,r}(0,\ldots,0,a,0,\ldots,0,b,0,\ldots,0)$, where $a$ and $b$ denote the $p$th and
$q$th strips of the $r$-circulant matrix $C_{n,r}$, respectively.
\\Clearly, the $(1,m+1)$ entry of the semicirculant matrix $[C_{n,r}]^{k}$ is
\begin{equation}\label{eq 1112}
[C_{n,r}]^{k}_{m}=
a^{k-\frac{m-kp}{q-p}}b^{\frac{m-kp}{q-p}}\displaystyle\binom{k}{\frac{m-kp}{q-p}},
\end{equation}
where we have adopted the convention that if $s\not\in \mathbb{N}$, then $\binom{k}{s}=0$.
\\From Theorem~\ref{Thm C} and formula~\eqref{eq 1112}, the ith strep of $r$-circulant matrix $C_{n,r}^{k}$ is
\begin{equation}
(C_{n,r}^{k})_{i}=\sum_{m\equiv i\; (\bmod\, n)} a^{k-\frac{m-kp}{q-p}}b^{\frac{m-kp}{q-p}}\displaystyle
\binom{k}{\frac{m-kp}{q-p}}r^{\floor*{\frac{m}{n}}}.\end{equation}
Especially, for $p=1$ and $q=n-1$, we find the following formula that gives the ith strep of the kth power of
$r$-circulant matrix $C_{n,r}=\text{circ}_{n}(0,a,0,\ldots,0,b)$
\begin{equation}\label{eq 0100}(C_{n,r}^{k})_{i}=\sum_{m\equiv i\; (\bmod\, n)} a^{k-\frac{m-k}{n-2}}b^{\frac{m-k}{n-2}}
\displaystyle\binom{k}{\frac{m-k}{n-2}}r^{\floor*{\frac{m}{n}}},\end{equation}
which is the same as that found by Jiang~\cite{Jiang} with $h$ replaced by $\displaystyle\frac{m-k}{n-2}$.
\end{ex}
\section{Conclusion}
We proposed a method that does not require any assumption to compute the arbitrary positive integer powers for
semicirculant matrices over an arbitrary unitary commutative ring. Accordingly, we derived an easy method to
compute the arbitrary positive integer powers for $r$-circulant matrices. The advantage of this method is that it does not
require to solve the difficult problem of determining the solution set of equation~\ref{eq;;010}. By comparing this method with~\cite{Feng, Jiang, Dav, Guti1, Guti2, Kok, Rim1, Rim2, Rim3, Rim4, Zhao} we see clearly that our method is much simpler than those used in~\cite{Feng, Jiang, Dav, Guti1, Guti2, Kok, Rim1, Rim2, Rim3, Rim4, Zhao} .


\begin{thebibliography}{99}
\bibitem{Feng} J. Feng, \emph{A note on computing of positive integer powers for circulant matrices,} Appl. Math. Comput. 223 (2013) 472--475, \url{doi:http://dx.doi.org/10.1016/j.amc.2013.08.016}.

\bibitem{Jiang} Z. Jiang, H. Xin, H. Wang, \emph{On computing of positive integer powers for r-circulant matrices,} Appl. Math. Comput. 265(C) (2015)
409--413, \url{doi:http://dx.doi.org/10.1016/j.amc.2015.05.022}.

\bibitem{Dav} P. Davis, \emph{Circulant Matrices,} Ams Chelsea Publishing, Providence, Rhode Island, 2012.

\bibitem{Guti1} J. Guti\'{e}rrez-Guti\'{e}rrez, \emph{Positive integer powers of complex skew-symmetric circulant matrices,} Appl. Math. Comput. 202(2)
(2008) 798--802, \url{doi:10.1016/j.amc.2008.03.024}.

\bibitem{Guti2} J. Guti\'{e}rrez-Guti\'{e}rrez, \emph{Positive integer powers of complex symmetric circulant matrices,} Appl. Math. Comput. 202(2) (2008)
877--881, \url{doi:10.1016/j.amc.2008.02.010}.

\bibitem{Kok} F. K\"{o}ken, D. Bozkurt, \emph{Positive integer powers for one type of odd order circulant matrices,} Appl. Math. Comput. 217(9) (2011)
4377--4381, \url{doi:10.1016/j.amc.2010.10.030}.

\bibitem{Rim1} J. Rimas, \emph{On computing of arbitrary positive integer powers for one type of odd order symmetric circulant matrices-{I},}
Appl. Math. Comput. 165(1) (2005) 137--141, \url{doi:10.1016/j.amc.2004.04.023}.

\bibitem{Rim2} J. Rimas, \emph{On computing of arbitrary positive integer powers for one type of odd order symmetric circulant matrices-{II},}
Appl. Math. Comput. 169(2) (2005) 1016--1027, \url{doi:10.1016/j.amc.2004.11.003}.

\bibitem{Rim3} J. Rimas, \emph{On computing of arbitrary positive integer powers for one type of even order symmetric circulant matrices-{I},}
Appl. Math. Comput. 172(1) (2006) 86--90, \url{doi:10.1016/j.amc.2005.01.121}.

\bibitem{Rim4} J. Rimas, \emph{On computing of arbitrary positive integer powers for one type of even order symmetric circulant matrices-{II},}
Appl. Math. Comput. 174(1) (2006) 511--523, \url{doi:10.1016/j.amc.2005.04.102}.

\bibitem{Zhao} G. Zhao, \emph{A cogredient algorithm for the m-th power of r-circulant matrices,} volume 2, International Conference on Computer Technology and Development (2009) 557--560. IEEE, \url{doi: 10.1109/ICCTD.2009.220}.

\bibitem{Hen} P. Henrici, \emph{Applied and Computational Complex Analysis,} volume 1, John Wiley, New York, 1974.

\bibitem{Brow} W. Brown, \emph{Matrices Over Commutative Rings,} Marcel Dekker, Inc., New York, 1993.

\bibitem{Wiki} Wikipedia, \emph{Formal power series,} \url{https://en.wikipedia.org/wiki/Formal_power_series}, 2020 (accessed May 31, 2020).
\end{thebibliography}
\end{document}